\documentclass[a4,11pt]{article}
\usepackage[centertags]{amsmath}
\usepackage{amsmath,amsfonts,amssymb,amsthm}
\usepackage[latin1]{inputenc}
\newcommand{\bb}{\begin{equation}}
\newcommand{\ee}{\end{equation}}

\newcommand{\la}{\label}

\newcommand{\lh}{\Delta_{H^{1}}}
\def\p{\partial}

\def\x1{{\xi }_{xx}}
\def\x2{{\xi }_{yy}}
\def\x3{{\xi }_{xy}}

\def\e1{{\eta }_{xx}}
\def\e2{{\eta }_{yy}}
\def\e3{{\eta }_{xy}}

\def\var{\varphi}

\def\la{\lambda}

\def\l1{{\lambda}_1}

\def\kd{\partial}

\begin{document}
\pagenumbering{arabic}
\title{\huge \bf Divergence Symmetries of Critical Kohn-\\Laplace Equations
on Heisenberg groups}
\author{\rm \large Yuri Bozhkov and Igor Leite Freire\\ \\
\it Instituto de Matem\'atica,
Estat\'\i stica e \\ \it Computa\c c\~ao Cient\'\i fica - IMECC \\
\it Universidade Estadual de Campinas - UNICAMP \\ \it C.P.
$6065$, $13083$-$970$ - Campinas - SP, Brasil
\\ \rm E-mail: bozhkov@ime.unicamp.br \\ \ \ \ \ \
igor@ime.unicamp.br }
\date{\ }
\maketitle
\vspace{1cm}
\begin{abstract}
We show that any Lie point symmetry of semilinear Kohn-Laplace
equations on the Heisenberg group $H^1$ with power nonlinearity is
a divergence symmetry if and only if the corresponding exponent
assumes critical value.
\end{abstract}

\

\section{Introduction}

\

 Recently a great renewed and increasing interest in the Heisenberg group has been
 manifested from
 both analysts and geometers. In this regard, in the last few decades  partial differential equations on the Heisenberg
 group were studied by various authors using different methods.
 There is a big variety of
 results among which we shall mention a few. In \cite{gl3} Garofalo and Lanconelli
 obtained existence, regularity and nonexistence results for
 semilinear PDEs involving Kohn-Laplace operators. General nonexistence results for
 solutions of differential inequalities on the Heisenberg group were established
 by Pokhozhaev and Veron in \cite{pv}. The work \cite{la} by Lanconelli is
 a survey on a series of results concerning
 critical semilinear equations on the Heisenberg
 group. For further details we direct the interested
 reader to these works and the references therein as well as to the existing
 internet instruments for search of mathematical information.

 In this paper we apply the S. Lie symmetry theory of
 differential equations (\cite{ ol,ov}) to the study of a model Kohn-Laplace equation
 on the Heisenberg group. Namely, we study the variational and divergence symmetries
 of the following differential equation on
 the Heisenberg group $H^1$:
 \bb\label{u01} \lh u +u^p=0,
 \ee
 where $\lh =X^2+Y^2 $ is the Kohn-Laplace operator,
 \[ X=\frac{\p}{\p x}+2y \frac{\p}{\p t} \]
 and
 \[ Y=\frac{\p}{\p y}-2x \frac{\p}{\p t} \]
 are the generators of \it left \rm multiplication in $H^1$. In more details, the equation (\ref{u01}) for
 $u=u(x,y,t): {\mathbb{R} }^3\rightarrow \mathbb{R} $ reads
 \bb\label{u02} u_{xx}+u_{yy}+4(x^{2}+y^{2})u_{tt}+4yu_{xt}-4xu_{yt}+u^p=0 .\ee

 In a previous work \cite{yi1} we obtained a complete group classification of semilinear Kohn-Laplace
 equations on $H^1$. In the case of nonlinearity of power type, the result states that the symmetry group of
 (\ref{u01}) for $p\neq 0$, $p\neq 1$ consists of translations in $t$, rotations in the $x-y$ plane,
 right multiplications in the Heisenberg group $H^1$ and a dilation, generated respectively by
  \bb\label{u03} T=\frac{\kd }{\kd t},\;\;\; R= y\frac{\kd }{\kd
 x}-x \frac{\kd }{\kd y},\;\;\; \tilde{X}=\frac{\p}{\p x}-2y \frac{\p}{\p
 t}, \;\;\;\tilde{Y}=\frac{\p}{\p y}+2x \frac{\p}{\p t} \ee
 and
 \bb\label{u04} Z= x\frac{\p}{\p x}+y\frac{\p}{\p y}+2t
 \frac{\p}{\p t}+\frac{2}{1-p}u\frac{\p}{\p u}. \ee
 Moreover, if $p=3$ the symmetry group can be expanded by the following generators:
 \bb
 \label{u05} V_1 = (xt-x^{2}y-y^{3})\frac{\p }{\p x} +
 (yt+x^{3}+xy^{2})\frac{\p }{\p y} + (t^{2}-(x^{2}+y^{2})^{2})\frac{\p }{\p t}-
 t u \frac{\p }{\p u} ,\ee
  \bb\label{u06} V_2  = (t-4xy)\frac{\p }{\p x} +
 (3x^{2}-y^{2})\frac{\p }{\p y} - (2yt+2x^{3}+2xy^{2})\frac{\p }{\p t}
 + 2 y u \frac{\p }{\p u} ,\ee
  \bb \label{u07} V_3 = (x^{2}-3y^{2})\frac{\p }{\p x} +
 (t+4xy)\frac{\p }{\p y} + (2xt-2x^{2}y-2y^{3})\frac{\p }{\p t}
 - 2x u \frac{\p }{\p u} .\ee

 The purpose of this paper is to find out which of the above symmetries are variational or divergence
 symmetries.

 We denote by $G$ the five-parameter Lie group of point transformations generated by $T,R, \tilde{X},
 \tilde{Y}$ and $Z$. Then  our first result can be formulated as follows.

 \

 \bf Theorem 1. \it The Lie point symmetry group $G$ of the Kohn-Laplace equation $(\ref{u01})$ is a variational
 symmetry group if and only if $p=3$.

 \

 \rm We recall that the homogeneous dimension of the Heisenberg group $H^n$ is given by $Q=2n+2$ and that the
 critical Sobolev exponent is $(Q+2)/(Q-2)$. Hence the Theorem 1 means that the dilation $Z$ is a
 variational symmetry if and only if $p$ equals to the critical exponent. Actually, the latter property
 holds for $H^n$, $n>1$, and we shall come back to this point later.

 Further we show that in the critical case the additional symmetries $V_1,V_2,V_3$ are divergence
 symmetries. For this purpose we find explicitly the vector-valued `potentials' which determine $V_1,V_2$
 and $V_3$ as divergence symmetries.

 The main result in this paper is the following

 \

 \bf Theorem 2. \it Any Lie point symmetry of the Kohn-Laplace equation
 \bb \label{u08} \lh u + u^3 =0 \ee
 is a divergence symmetry.

 \

 \rm As it is well-known, the divergence symmetries determine conservation laws via the
 Noether Theorem \cite{ol}. Thus the next steps in this research would be to
 establish the conservation laws corresponding
 to the already studied variational and divergence symmetries as well as to study the invariant solutions of
 the Kohn-Laplace equations. These problems will be treated elsewhere \cite{yi2,yi3}.

 We observe that, by Theorem 2, all Lie point symmetries of the Kohn-Laplace equation (\ref{u08}) are
 divergence symmetries. This confirms the validity of the general property, established and discussed in
 \cite{yb1}, stating that the Lie point symmetries of critical quasilinear differential equations with power
 nonlinearities are divergence symmetries. The fact that this property should be valid for differential
 equations on Heisenberg groups was conjectured by Enzo Mitidieri in June 2003 \cite{em}.

 This paper is organized as follows. In sections 2 and 3 we prove theorems 1 and 2 respectively. In section
 3 we discuss the generalization of the obtained results to the Heisenberg group $H^n$,
 $n>1$.

 \

 \section{The variational symmetries}

 \

 To begin with, we note that the Kohn-Laplace equation
 \bb\label{u09} \lh u + f(u) =0 \ee
 is the Euler-Lagrange equation for the functional
 \[ \int L(x,y,t,u,u_x, u_y,u_t)dxdydt, \]
 where the integration is performed over ${\mathbb{R} }^3$, the function of Lagrange is given by
 \[ L=\frac{1}{2}(Xu)^2 +\frac{1}{2}(Yu)^2 -\int_{0}^{u} f(s)ds,\]
 or, equivalently,
 \bb\label{u10} L=\frac{1}{2} u_x^2 +\frac{1}{2} u_y^2 + 2(x^2+y^2)u_t^2 +2yu_xu_t-2xu_yu_t -
 \int_{0}^{u} f(s)ds, \ee
 and the function $u$ is assumed to satisfy appropriate vanishing conditions as
 $d=(t^2+(x^2+y^2))^{1/4}\rightarrow\infty $.

 \bf Proof of Theorem 1. \rm By the general theory of symmetries of differential equations
 \cite{ol} it is enough to show that the generators of $G$ determine variational symmetries.
 Indeed, by the infinitesimal criterion for invariance \cite{ol}, p. 257, $G$ is a variational
 symmetry group if and only if
 \bb\label{u11} W^{(1)} L + L (D_x \xi + D_y\phi +D_t\tau )=0 \ee
 for all $ (x,y,t,u,u_x, u_y,u_t)$ and every infinitesimal generator
 \[ W=\xi \frac{\p }{\p x}+\phi \frac{\p }{\p y} +\tau \frac{\p }{\p t} +
 \eta \frac{\p }{\p u}.\]
 (Recall that $W^{(1)}$ is the first order extension of $W$, see \cite{ol}.)

 Aiming to verify (\ref{u11}) for $T,R,\tilde{X}, \tilde{Y}$ and $Z$, we first calculate the
 corresponding first order extensions using the formulae for the extended infinitesimals \cite{ol}:
 \bb\label{u12} T^{(1)}=T, \ee
 \bb\label{u13} R^{(1)}=R+u_y\frac{\p }{\p u_x}- u_x\frac{\p }{\p u_y}, \ee
 \bb\label{u14} {\tilde{X}}^{(1)}=\tilde{X} + 2u_t\frac{\p }{\p u_y}, \ee
 \bb\label{u15} {\tilde{Y}}^{(1)}=\tilde{Y} - 2u_t\frac{\p }{\p u_x}, \ee
 and
 \bb\label{u16} Z^{(1)} =Z + \frac{1+p}{1-p}u_x\frac{\p }{\p u_x} +
 \frac{1+p}{1-p}u_y\frac{\p }{\p u_y} + \frac{2p}{1-p}u_t\frac{\p }{\p u_t} .\ee

 Then from  (\ref{u03}), (\ref{u10}), (\ref{u12}), (\ref{u13}), (\ref{u14}), (\ref{u15})
 we obtain easily that $T,R,\tilde{X}, \tilde{Y}$  satisfy (\ref{u11}). Hence they determine
 variational symmetries for arbitrary $f(u)$.

 Further, let $\xi =x, \phi =y, \tau = 2t, \eta =2u/(1-p)$ be the infinitesimals of the dilation $Z$.
 Then the left-hand side of (\ref{u11}) with $W=Z$ reads
 \[ Z^{(1)}L+L (D_x \xi +D_y\phi +D_t\tau )=\left [
 x\frac{\p}{\p x}+y\frac{\p}{\p y}+2t
 \frac{\p}{\p t}+\frac{2}{1-p}u\frac{\p}{\p u} \right. \]
 \[ \left. + \frac{1+p}{1-p}u_x\frac{\p }{\p u_x} +
 \frac{1+p}{1-p}u_y\frac{\p }{\p u_y} + \frac{2p}{1-p}u_t\frac{\p }{\p u_t} \right ]L +4 L\]
 where
 \[ L=\frac{1}{2} u_x^2 +\frac{1}{2} u_y^2 + 2(x^2+y^2)u_t^2 +2yu_xu_t-2xu_yu_t -
 \frac{1}{p+1} u^{p+1}.\]
 After a differentiation and simplifying we obtain
 \[ Z^{(1)}L+L (D_x \xi +D_y\phi +D_t\tau )= \]
 \[ \frac{3-p}{1-p}(u_x^2+u_y^2+4(x^2+y^2)u_t^2+2yu_xu_t-2xu_yu_t) +\frac{2(3-p)}{p^2-1}u^{p+1}.\]
 Hence $Z$ is a variational symmetry if and only if $p=3$, which completes the proof of
 Theorem 1.

 \

 \section{The divergence symmetries}

 \

 In this section we prove Theorem 2.

 Recall that a point transformation with infinitesimal generator
 \[ W=\xi \frac{\p }{\p x}+\phi \frac{\p }{\p y} +\tau \frac{\p }{\p t} +
 \eta \frac{\p }{\p u} \]
 is a divergence symmetry for $\int L$ if and only if there exists a vector function
 $\var =({\var }_1 , {\var }_2, {\var }_3)$ of $u$ and its derivatives up to some finite order such
 that
 \bb\label{u17} W^{(1)} L + L (D_x \xi + D_y\phi +D_t\tau )= Div\;(\var ). \ee
 Since the variational symmetries are divergence symmetries with $\var =0 $, by Theorem 1 it is enough
 to prove that $V_1, V_2,V_3$ are divergence symmetries. For this purpose we shall find the
 corresponding `potentials' $\var $.

 For the symmetry $V_1$ we have
 \[ \xi = xt-x^{2}y-y^{3}, \;\;\; \phi = yt+x^{3}+xy^{2}, \]
 \[ \tau = t^{2}-(x^{2}+y^{2})^{2}, \;\;\; \eta = -tu .\]
 We calculate the first order extension of $V_1$:
 \[ V_1^{(1)}= V_1+ {\eta }^{(1)}_x \frac{\p }{\p u_x} +
  {\eta }^{(1)}_y \frac{\p }{\p u_y} +{\eta }^{(1)}_t \frac{\p }{\p u_t} ,\]
 where the extended infinitesimals are given by
 \[ {\eta }^{(1)}_x= 2(xy-t)u_x-(3x^2+y^2)u_y+4x(x^2+y^2)u_t,\]
 \[ {\eta }^{(1)}_y= (x^2+3y^2)u_x-2(t+xy)u_y+4y(x^2+y^2)u_t,\]
 \[ {\eta }^{(1)}_t= -u - xu_x-yu_y-3tu_t.\]
 Then, after some tedious work, we obtain
 \bb\label{v17} V_1^{(1)}L + L (D_x \xi +D_y\phi +D_t\tau )= 2xuu_y-2yuu_x-4(x^2+y^2)u_t . \ee
 Thus $V_1$ is not a variational symmetry. Let
 \bb\label{u18} A_1=-yu^2,\;\;\; A_2=xu^2,\;\;\; A_3=-2(x^2+y^2)u^2. \ee
 Then from (\ref{v17}) and (\ref{u18})
 \[ V_1^{(1)}L + L (D_x \xi +D_y\phi +D_t\tau )=Div(A), \]
 where $A=(A_1,A_2,A_3)$. Hence $V_1$ is a divergence symmetry.

 Analogously, for the symmetries $V_2$ and $V_3$, after another tedious work, we obtain
 \bb\label{u19} V_2^{(1)}L + L (D_x \xi +D_y\phi +D_t\tau )= 2u u_y-4x u u_t\ee
 and
 \bb\label{u20} V_3^{(1)}L + L (D_x \xi +D_y\phi +D_t\tau )=-2u u_x-4y u u_t.\ee
 Now we denote
 \bb\label{u21} B=(0,u^2,-2x u^2)\ee
 and \bb\label{u22} C=(-u^2,0,-2y u^2) .\ee
 Then by (\ref{u19})-(\ref{u22}) we see that $V_2$ and $V_3$ satisfy (\ref{u17}) with $B$ and $C$,
 respectively, in the place of $\var $. Thus $V_2$ and $V_3$ are divergence symmetries.

 \

 \section{On the generalization to $H^n$, $n>1$ }

 \

 In this section we comment briefly on the generalization of our approach to the Heisenberg group
 $H^n$, $n>1$.

 First we observe that the Kohn-Laplace equation
 \bb\label{u23} {\Delta }_{H^n} u+f(u)=0, \ee
 or equivalently
 \[ u_{x_i x_i}+u_{y_i y_i}+4(x^{2}_i+y^{2}_i)u_{tt}+4y_i u_{x_i t}-4x_i u_{y_i t}+f(u)=0 \]
 is the Euler-Lagrange equation of the functional
 \[ J[u]=\int L,\]
 with
 \[ L=\frac{1}{2}(X_i u)^2 +\frac{1}{2}(Y_i u)^2 -\int_{0}^{u} f(s)ds\]
 \[ = \frac{1}{2} u_{x_i}^2 +\frac{1}{2} u_{y_i}^2 + 2(x^2_i+y^2_i)u_t^2 +2y_i u_{x_i}u_t-
 2x_i u_{y_i}u_t - \int_{0}^{u} f(s)ds,\]
 where $ X_i=\frac{\p}{\p x_i}+2y_i \frac{\p}{\p t} $, $Y_i=\frac{\p}{\p y_i}-2x_i \frac{\p}{\p t} $
 and a summation over $i=1,2,...n$ is assumed.

 Using the definition of Lie point symmetry of a differential equation, one can show that the scaling
 transformation
 \[ \left\{\begin{array}{lll} x_j^{*} &=&\la x_j, \\
 y_j^{*} & = &\la y_j, \\
 t^{*} &=&{\la }^2t, \\
 u^{*}&=&{\la }^{\frac{2}{1-p}} u \end{array} \right. \]
 is admitted by the equation
 \bb\label{u24} {\Delta }_{H^n} u+ u^p=0. \ee

 Then performing this change of variables in the functional $J$ it is easy to see that the dilation
 \[ Z= x_i \frac{\p}{\p x_i}+y_i \frac{\p}{\p y_i}+2t
 \frac{\p}{\p t}+\frac{2}{1-p}u\frac{\p}{\p u} \]
 is a variational symmetry if and only if
 \[ p=\frac{n+2}{n}=\frac{Q+2}{Q-2}. \]
 Thus the equation (\ref{u24}) admits a variational symmetry group containing $Z$ if and only if $p$
 assumes the critical value.

 In conclusion, we note that the complete group classification of the Kohn-Laplace equations
 (\ref{u23}) is available only for $n=1$ (\cite{yi1}). We conjecture that in the critical case
 the symmetry group of (\ref{u24}) can be expanded and, in addition, that all Lie point
 symmetries in
 this case are divergence symmetries.
 \

 \begin{center}\textbf{Acknowledgements}\end{center}

 We express our deep gratitude to Professor Enzo Mitidieri for useful discussions
 and encouragement. Y. Bozhkov would also like to thank FAPESP and CNPq,
 Brasil, for financial support. I. L. Freire is grateful to
 CAPES, Brasil, for financial support.

\end{document}